\documentclass[12pt]{article}
\usepackage{amsfonts,hyperref,amsmath}
\usepackage[latin1]{inputenc}
\usepackage[english]{babel}
\usepackage{graphicx}
\usepackage{spverbatim}
\usepackage{hyperref}
\usepackage[top=20mm, bottom=20mm, left=20mm , right=20mm]{geometry}

\begin{document}

\renewcommand{\theequation}{\thesection.\arabic{equation}}
\newcommand{\newsection}{    
\setcounter{equation}{0}
\section}
\renewcommand{\appendix}[1]{
    \addtocounter{section}{1}
    \setcounter{equation}{0}
    \renewcommand{\thesection}{\Alph{section}}
    \section*{Appendix \thesection\protect\indent #1}
    \addcontentsline{toc}{section}{Appendix \thesection\ \ \ #1}
}
\newcommand\encadremath[1]{\vbox{\hrule\hbox{\vrule\kern8pt 
\vbox{\kern8pt \hbox{$\displaystyle #1$}\kern8pt} 
\kern8pt\vrule}\hrule}}
\def\enca#1{\vbox{\hrule\hbox{
\vrule\kern8pt\vbox{\kern8pt \hbox{$\displaystyle #1$}
\kern8pt} \kern8pt\vrule}\hrule}}

\newcommand\figureframex[3]{
\begin{figure}[bth]
\hrule\hbox{\vrule\kern8pt 
\vbox{\kern8pt \vbox{
\begin{center}
{\mbox{\epsfxsize=#1.truecm\epsfbox{#2}}}
\end{center}
\caption{#3}
}\kern8pt} 
\kern8pt\vrule}\hrule
\end{figure}
}
\newcommand\figureframey[3]{
\begin{figure}[bth]
\hrule\hbox{\vrule\kern8pt 
\vbox{\kern8pt \vbox{
\begin{center}
{\mbox{\epsfysize=#1.truecm\epsfbox{#2}}}
\end{center}
\caption{#3}
}\kern8pt} 
\kern8pt\vrule}\hrule
\end{figure}
}

\makeatletter
\@addtoreset{equation}{section}
\makeatother
\newtheorem{theorem}{Theorem}[section]
\newtheorem{conjecture}{Conjecture}[section]
\newtheorem{remark}{Remark}[section]
\newtheorem{proposition}{Proposition}[section]
\newtheorem{lemma}{Lemma}[section]
\newtheorem{corollary}{Corollary}[section]
\newtheorem{definition}{Definition}[section]
\def\br{\begin{remark}\rm\small}
\def\er{\end{remark}}
\def\bt{\begin{theorem}}
\def\et{\end{theorem}}
\def\bd{\begin{definition}}
\def\ed{\end{definition}}
\def\bp{\begin{proposition}}
\def\ep{\end{proposition}}
\def\bl{\begin{lemma}}
\def\el{\end{lemma}}
\def\bc{\begin{corollary}}
\def\ec{\end{corollary}}
\def\beaq{\begin{eqnarray}}
\def\eeaq{\end{eqnarray}}
\newcommand{\proof}[1]{{\noindent \bf proof:}\par
{#1} $\square$}

\newcommand{\rf}[1]{(\ref{#1})}
\newcommand{\eq}[1]{Eq.~(\ref{#1})}
\newcommand{\rfig}[1]{fig.~\ref{#1}}

\newcommand{\equ}[2]{\begin{equation}{\label{#1}}{#2}\end{equation}}

\newcommand{\beq}{\begin{equation}}
\newcommand{\eeq}{\end{equation}}
\newcommand{\beqq}{\begin{equation*}}
\newcommand{\eeqq}{\end{equation*}}
\newcommand{\bea}{\begin{eqnarray}}
\newcommand{\eea}{\end{eqnarray}}
\newcommand{\beaa}{\begin{eqnarray*}}
\newcommand{\eeaa}{\end{eqnarray*}}

\newcommand\eol{\hspace*{\fill}\linebreak}
\newcommand\eop{\vspace*{\fill}\pagebreak}

\newcommand{\hs}{\hspace{0.7cm}}
\newcommand{\vs}{\vspace{0.7cm}}
%
%

%
\renewcommand{\and}{{\qquad {\rm and} \qquad}}
\newcommand{\where}{{\qquad {\rm where} \qquad}}
\newcommand{\with}{{\qquad {\rm with} \qquad}}
\newcommand{\for}{{\qquad {\rm for} \qquad}}
\newcommand{\virg}{{\qquad , \qquad}}

\newcommand{\ie}{{\it i.e.}\ }


\newcommand{\Det}{{\,\rm Det}} \newcommand{\Tr}{{\,\rm Tr}\:}
\newcommand{\tr}{{\,\rm tr}\:}
\newcommand{\cte}{{\,\rm cte}\,}
\newcommand{\Res}{\mathop{\,\rm Res\,}}

\newcommand{\td}[1]{{\tilde{#1}}}

\renewcommand{\l}{\lambda}
\newcommand{\om}{\omega}
\newcommand{\calP}{{\cal P}}

\newcommand{\ii}{{\mathrm{i}}}
\newcommand{\e}{{\,\rm e}\,}
\newcommand{\ee}[1]{{{\rm e}^{#1}}}

\renewcommand{\d}{{{\partial}}}
\newcommand{\D}{{{\hbox{d}}}}
\newcommand{\dmat}[2]{\mathrm{d}_{\scriptscriptstyle{#1}}[#2]}

\newcommand{\Pint}{{\int\kern -1.em -\kern-.25em}} 
\newcommand{\Vol}{\mathrm{Vol}}

\newcommand{\moy}[1]{\left<{#1}\right>}

\renewcommand{\Re}{{\mathrm{Re}}}
\renewcommand{\Im}{{\mathrm{Im}}}

\newcommand{\sn}{{\rm sn}}
\newcommand{\cn}{{\rm cn}}
\newcommand{\dn}{{\rm dn}}
\newcommand{\ssq}[1]{{\sqrt{\sigma({#1})}}}

\newcommand{\hd}{{\hat\delta}}

\renewcommand{\l}{\lambda}
\renewcommand{\L}{\Lambda}
\renewcommand{\ssq}[1]{{\sqrt{\sigma({#1})}}}
\newcommand{\ovl}{\overline}

\begin{center}
\vspace{5pt}
{\large \bf {Locks and keys: How fast can you open several locks with too many keys?}}
\end{center}
\vspace{2pt}
\begin{center}\begin{tabular}{c}
\textbf{O. Marchal}\\
\textit{Universit\'e de Lyon, CNRS UMR 5208,} \\
\textit{Universit\'e Jean Monnet, Institut Camille Jordan, France \footnote{olivier.marchal@univ-st-etienne.fr}}
\end{tabular}
\end{center}

\vspace{20pt}

{\bf Abstract}: This short note is the result of a French ``Hippocampe internship''\footnote{``Hippocampe'' (sea-horse or hippocampus in French) internships are three days sessions during which small groups of students specialized in mathematics and physics are presented with a problem that they have to solve under the supervision of an academic professor.}
 that aims at introducing the world of research to young undergraduate French students. The problem studied is the following: imagine yourself locked in a cage barred with $n$ different locks. You are given a keyring with $N \geq n$ keys containing the $n$ keys that open the locks. In average, how many trials are required to open all locks and get out? The article studies $3$ different strategies and compare them. Implementation of the strategies are also proposed as illustrations of the theoretical results.

\section{Statement of the problem}
The precise statement of the problem and the various assumptions are the following:
\medskip

\underline{Statement of the problem}: \textit{We need to open $n$ locks with a keyring containing $N\geq n$ keys. We are interested in determining the average and the standard deviation of the number of trials to open all locks. The special case when the number of keys equals the number of locks ($n=N$) is included.}

\medskip

\underline{Assumptions}: The following assumptions about locks and keys are made:
\begin{itemize}\item All locks and keys are assumed to be indistinguishable.
\item All locks are different and require a specific key to be opened.
\item A lock can always be opened, i.e. the key required to open it belongs to the keyring.
\item A key cannot open more than one lock.
\end{itemize}

Consequently, the keyring contains the $n$ different keys corresponding to the locks as well as additional keys that do not open anything when $N>n$. The keys are randomly placed in the keyring. We denote by $T_{n,N}$ the number of trials to open all locks. We are interested in gathering as much information as possible on $T_{n,N}$ especially $E(T_{n,N})$ and $\text{Var}(T_{n,N})$.

\section{Strategy $1$: The totally random strategy} 
The first strategy one can think of is to select randomly one of the remaining keys and try it randomly into one the remaining locks to open. In particular we do not keep track about the previous tries except for the fact that we discard locks that are already opened and their corresponding keys. Let us introduce $X_i$ ($1 \leq i\leq n$) the number of trials required between the opening of $i-1^{\text{th}}$ lock and the $i^{\text{th}}$ lock. (for $i=1$, $X_1$ is just the number of trials required to open the first lock). Thus, we have:
\beq \label{T} T_{n,N}=\sum_{i=1}^n X_i\eeq
In this strategy, the variables $X_i$ are obviously independent since we do not keep track of the history of the trials. To open the $i^{\text{th}}$ lock, we choose randomly into the $N-i+1$ remaining keys. Therefore we get:
\beq \forall \,k\geq 1\,:\,  P(X_i=k)=\left(\prod_{i=1}^{k-1}\frac{N-i}{N-i+1}\right) \frac{1}{N-i+1}=\left(\frac{N-i}{N-i+1}\right)^{k-1}\frac{1}{N-i+1}\eeq
Thus we observe that $X_i$ is distributed as a geometric variable with parameters $q_i=\frac{N-i}{N-i+1}$. This is of course coherent with the fact that the strategy is memory-less. It is then straightforward to obtain the mean and standard deviation of $X_i$:
\beq E(X_i)=N-i+1 \text{  and  } \text{Var}(X_i)=(N-i)(N-i+1)\eeq
Since the variables $X_i$ are independent, we get from \eqref{T} that:
\footnotesize{\beq E(T_{n,N})=\sum_{i=1}^{n} E(X_i)=\frac{n(2N-n+1)}{2} \text{  and  } \text{Var}(T_{n,N})=\sum_{i=1}^{n} \text{Var}(X_i)=\frac{n(3N^2-3Nn+n^2-1)}{3}\eeq}\normalsize{}
In the case when $N=n$, we get:
\beq E(T_{n,n})=\frac{n(n+1)}{2} \text{ and } \text{Var}(T_{n,n})=\frac{n(n-1)(n+1)}{3}\eeq
Simulations (See \ref{Simulations}) show that $T_{n,N}$ is very well approached by a $\Gamma$ distribution with parameters $\Gamma\left(k=\frac{3n(2N-n+1)^2}{4(3N^2-3Nn+n^2-1)},\theta=\frac{2(3N^2-3Nn+n^2-1)}{3(2N-n+1)}\right)$. Parameters have been adjusted here so that expectation and standard deviation match (we recall that $E(\Gamma(k,\theta))=k\theta$ et $Var(\Gamma(k,\theta))=k\theta^2$). Appearance of the $\Gamma$ distribution is not surprising here. Indeed, geometric variables are a discrete version of the exponential distribution (typical of memory-less phenomenon) that is a special case of $\Gamma$ distributions $(k=1)$. Additionally, the sum of independent $\Gamma$ distributions with parameters $(k_i,\theta)$ is known to remain a $\Gamma$ distribution with parameters $\left(k=\underset{i=1}{\overset{n}{\sum}} k_i, \theta\right)$. However note that our parameters $\theta_i$ are different so that the last theorem cannot be applied properly but only in an approximate manner.

\medskip

Of course, the totally random strategy is not the optimal strategy since we do not keep track of the combinations tried previously. These strategies are developed in the next sections.

\section{Strategy $2$: Take a lock and test all keys\label{Stra1}}
\underline{Strategy}: The strategy is now the following: we choose a lock and then we try all remaining (that is to say all keys except those that have already opened a lock) keys successively until we find the (unique) right key that opens the lock. Then we continue with the next lock by starting over with the first remaining key and so on until all are opened.

\underline{Remark}: Since a key opens at most one lock, we can discard the successful keys. Moreover when a lock is opened, all remaining keys (including those that have been tried and those that have not been tried on the lock) remain equivalent for the next locks. Indeed, since we now that a lock is exactly opened by one key, we get that when a lock is opened, all keys that have not been tried immediately acquire the property that they do not open this lock, a property that was already known for the keys that have been tried on the lock. In particular the way we will order the keys for the next lock does not affect the distribution of the number of trials.

\medskip

As in the previous section, let $X_i$ ($1 \leq i\leq n$) be the number of trials required between the opening of $i-1^{\text{th}}$ lock and the $i^{\text{th}}$ lock. (for $i=1$, $X_1$ is just the number of trials required to open the first lock). Equation \eqref{T} still holds. Moreover from the previous remark, we get that the variables $X_i$ are still independent. To open the $i^{\text{th}}$ lock, only $N+1-i$ keys remain and we get for $1\leq k\leq N+1-i$:
\beq P(X_i=k)=\frac{N-i}{N+1-i}\frac{N-i-1}{N-i}\dots \frac{N+1-i-k}{N+2-i-k}\frac{1}{N+1-i-k}=\frac{1}{N+1-i}\eeq
Consequently we observe that $X_i$ follows a discrete uniform distribution:
\beq X_i \sim \text{Uniform}\left(\{1,\dots,N-i+1\}\right)\eeq
In the end we get:
\beq E(X_i)=\frac{N+2-i}{2} \text{  and  } \text{Var}(X_i)=\frac{(N+2-i)(N-i)}{12}\eeq
and eventually
\beq \label{Res} E(T_{n,N})=\frac{n(2N+3-n)}{4}\,\, \text{ and } \text{Var}(T_{n,N})=\frac{n(6N^2-6Nn+6N+2n^2-3n-5)}{72}\eeq
In the case when $N=n$ we get:
\beq \label{Res2} E(T_{n,n})=\frac{n(n+3)}{4}\,\, \text{ and } \text{Var}(T_{n,n})=\frac{n(2n+5)(n-1)}{72}\eeq

Simple computations shows that this strategy is in average better than the totally-random one even in the case when $N=n$. Simulations (See section \ref{Simulations}) of this strategy are also in agreement with the theoretical results.

\section{Strategy $3$: Take a key and try it in all locks \label{Stra2}}
\underline{Strategy}: We now consider the following strategy: we take a key and we try it into all remaining locks. We stop when it opens a lock or when we have exhausted all remaining locks. We then take the next key and start over with the first remaining lock and so on until all locks are opened.

\medskip

Let us denote as before $X_i$ ($1 \leq i\leq n$) the number of trials required between the opening of $i-1^{\text{th}}$ lock and the $i^{\text{th}}$ lock. (for $i=1$, $X_1$ is just the number of trials required to open the first lock) so that equation \eqref{T} still holds. In this situation, the variables $X_i$ are no longer independent. Indeed, let us assume that to open the first lock, we have used all $N-n$ useless keys (keys that do not open any locks), then it is now impossible that $X_2$ is greater than $n-1$ since we only have $n-1$ keys and $n-1$ locks remaining. On the contrary, if the first key opens a lock (that is $X_1\leq n$), then $X_2$ may be strictly greater than $n$ if we choose a useless key. Hence we see that the range of $X_2$ depends on the result of $X_1$, thus $X_1$ and $X_2$ cannot be independent. The case $N=n$ appears to be singular, since in that case the variables $X_i$ remain independent (since a key must open a lock). Determining the distribution of $X_i$ is more challenging too. Indeed $X_p=k$ depends on the number of keys that have been discarded until the opening of the $p-1^{\text{th}}$ lock. More specifically, if $l$ locks have been opened before the trial of the $p^{\text{th}}$ key, then the key opening lock $p$ must be in position $k+l$ in the keyring. Consequently we get:
\beq P(X_p=k)=\frac{1}{N}\sum_{l=0}^{p-1} \frac{ \binom{k+l-1}{l}\binom{N-k-l}{p-1-l}}{\binom{N-1}{p-1}}\eeq   
We want to prove the following:
\begin{proposition} \label{Prop}The following identity holds:
\beq \frac{1}{N}\sum_{l=0}^{p-1} \frac{ \binom{k+l-1}{l}\binom{N-k-l}{p-1-l}}{\binom{N-1}{p-1}}=\frac{1}{N-p+1} \,\,,\,\,\forall\, k\leq N-p+1\eeq
\end{proposition}

\medskip
\underline{Proof}:
This identity is a consequence of the Chu-Vandermonde identity\footnote{The Chu-Vandermonde identity can be proved easily by identification of the order $x^{N+1}$ in the binomial expansion of the identity $(1+x)^{N+m}=(1+x)^N(1+x)^m$}:
\beq \forall\, J\leq K\leq N \,:\,\, \binom{N+1}{K+1}=\sum_{m=0}^{K}\binom{m}{J}\binom{N-m}{K-J} \eeq
Then replacing $N$ by $N-1$ and using the symmetry of the binomial coefficients we get to:
\beq \forall \, J\leq K\leq N-1\,\,:\,\, \binom{N}{N-K-1}=\sum_{m=0}^{K}\binom{m}{J}\binom{N-1-m}{K-J} \eeq
Taking $K=N-p$ and $J=k-1$ we find
\beq  \forall \, k\leq N-p+1\,\,:\,\,\binom{N}{p-1}=\sum_{m=k-1}^{p+k-2}\binom{m}{k-1}\binom{N-1-m}{N-p-k+1}\eeq
where we restrict the sum so that only non-zero binomial coefficients are present. Performing the change $l=m-k+1$ gives: 
\beq \sum_{l=0}^{p-1}\binom{k+l-1}{k-1}\binom{N-k-l}{N-k-p+1}=\binom{N}{p-1}=\frac{N}{N-p+1}\binom{N-1}{p-1} \,\,,\,\,\forall\, k\leq N-p+1\eeq
This last equality is then equivalent to proposition \ref{Prop}.

\medskip

From proposition \ref{Prop} we conclude that the variable $X_p$ is distributed like a discrete uniform distribution:
\beq X_i\sim \text{Uniform}\left(\{1,\dots,N+1-i\}\right)\eeq
Expectation and standard deviation are then straightforward to compute:
\beq E(X_i)=\frac{(N+2-i)}{2} \text{ and }  \text{Var}(X_i)=\frac{(N+2-i)(N-i)}{12}\eeq
From \eqref{T} we finally get:
\beq E(T_{n,N})=\sum_{i=1}^n\frac{(N+2-i)}{2}=\frac{n(2N-n+3)}{4}\eeq
We recover the same result as in the previous strategy \eqref{Res}. In particular $T_{n,N}$ is still the sum of uniform variables. However in this approach, it is unclear to see why the variables $X_i$ should be uncorrelated. Consequently the computation of the variance is not directly possible except for the case $N=n$ where the variable $X_i$ are independent and \eqref{Res2} is recovered.

\medskip

The last result indicates that the last two strategies seem similar since their results are identical. However computations in the last strategy seem much more complicated than for the previous one and fails to explain why the strategy are equivalent. In fact, \textbf{the last two strategies are exactly equivalent no matter the order of keys selected: they will provide exactly the same number of trials for any initial conditions}. We will prove this result in the next section using a combinatorial one-to-one correspondence.

\section{Equivalence of strategies $2$ and $3$}
Let us order the locks from $1$ to $n$ and number the keys that open a lock with their corresponding lock number. The useless keys are then numbered randomly from $n+1$ to $N$ (the number itself is of no interest, the only piece of information that we need is that it does not open anything). Then an initial condition corresponds to the choice of a random permutation of $\{1,\dots, N\}$ of the keys. More specifically, the randomness of the problem consists in drawing a random permutation of $\{1,\dots,N\}$ uniformly on the symmetric group $S_N$ (with $N!$ possibilities). As soon as the permutation is given (i.e. the order of the keys given) then the number of trials is completely determined using one of the strategies developed in sections \ref{Stra1} or \ref{Stra2}. This defines properly the probability space on which we work on. Let us illustrate the situation on the following example:

\begin{center}\label{Fig1}
\includegraphics[width=16cm]{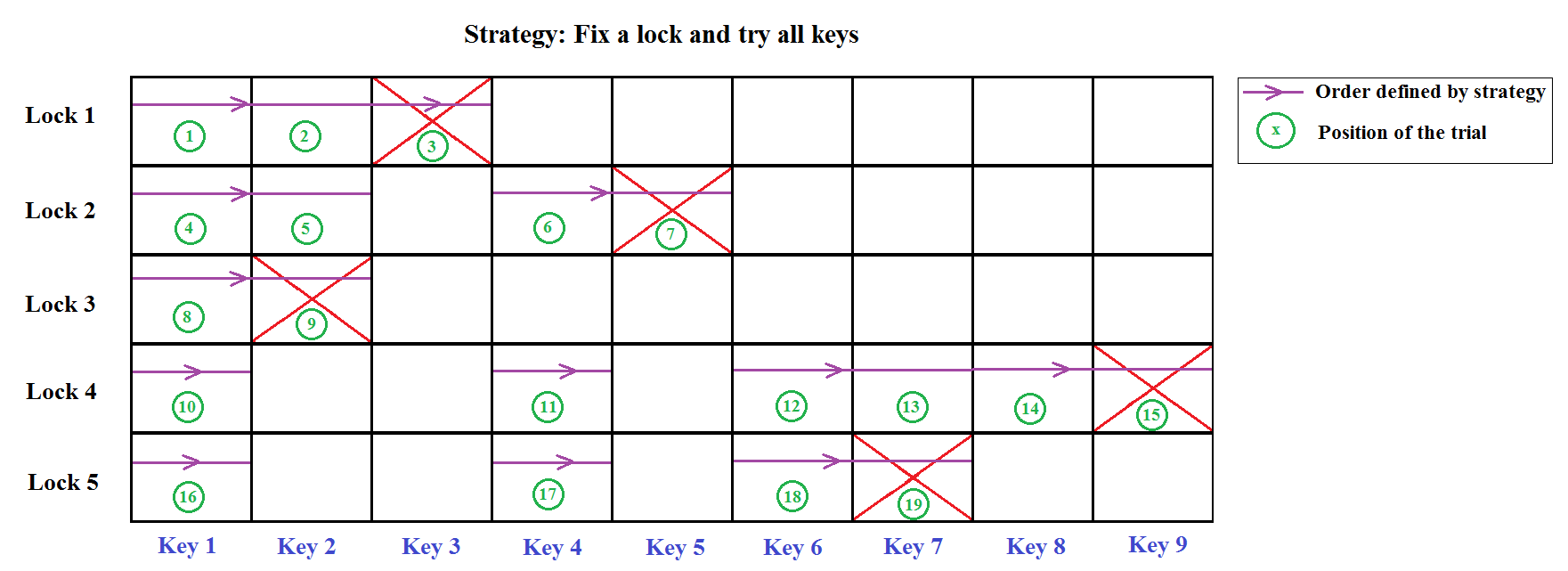}

\textit{Fig. $1$: Illustration of the strategy developed in section \ref{Stra1} when locks $1$ up to $5$ are opened respectively by keys in position $3,5,2,9,7$ in the keyring.}
\end{center}

\medskip

\begin{center}\label{Fig2}
\includegraphics[width=16cm]{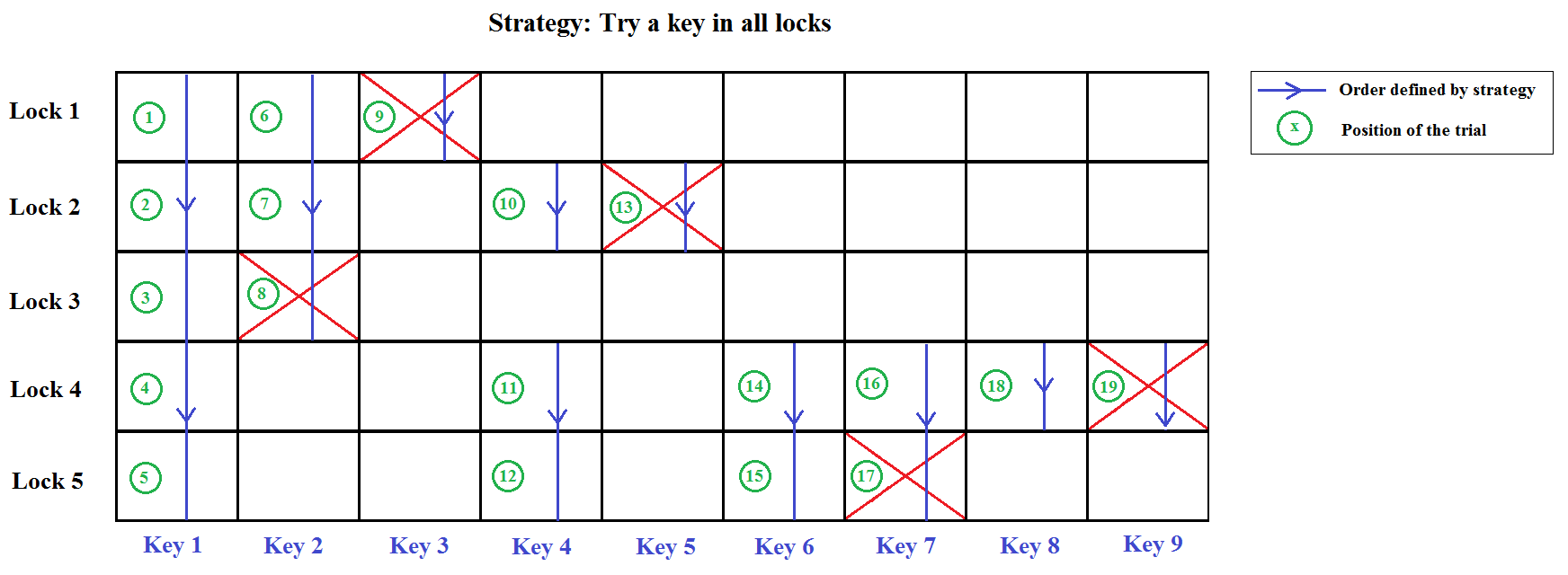}

\textit{Fig. $2$: Illustration of the strategy developed in section \ref{Stra2} when locks $1$ up to $5$ are opened respectively by keys in position $3,5,2,9,7$ in the keyring. Boxes ticked in red indicate the position of the keys in the keyring.}
\end{center}

\medskip
\medskip

For both strategies, it is easy to prove that the only boxes tried in the previous strategies are those being on the left and at the top of a ticked box. Indeed:
\begin{enumerate}\item In strategy developed in section \ref{Stra1}, boxes at the right of a ticked box will be discarded because we do not try the remaining keys as soon as a lock is open. Boxed below a ticked box will not be tried because as soon as a lock is opened, we do not test any other key on it.
\item In strategy developed in section \ref{Stra2}, boxes at the right of a ticked box will be discarded because as soon as a lock is opened we do not test any other key on it. Boxes below a ticked box will not be tried because as soon as a key has opened a lock, we do not try it on any other lock (since we know from the assumptions that it can at most open one lock).
\end{enumerate}

Hence figures $1$ and $2$ show that for any initial condition, the number of trials $T_{n,N}$ will be the same for both strategies: only the order of the various tries will differ in the two strategies.

\section{Numeric simulations \label{Simulations}}
Numeric simulations can be carried out to illustrate the various strategies. We present here the code in Maple as well as histograms of the distributions of $T_{n,N}$ for different values of $n$ and $N$. All strategies require a procedure to generate a random permutation:

\small{\begin{spverbatim}
RandomPerm:=proc(n)
local T:
T:=[seq(1..n)]:
return(Shuffle(T)):
end proc:
\end{spverbatim}} \normalsize{}
\medskip

Then, the totally random strategy can be implemented like:
\small{\begin{spverbatim}
RandomStrategy:=proc(NumberKeys,NumberLocks)
local NumberTrials, LocksRemaining, KeysRemaining, LockChosen, KeyChosen:
NumberTrials:=0:
LocksRemaining:=NumberLocks:
KeysRemaining:=NumberKeys:
while LocksRemaining>0 do
LockChosen:=rand(LocksRemaining)()+1:
KeyChosen:=rand(KeysRemaining)()+1:
if LockChosen=KeyChosen then LocksRemaining:=LocksRemaining-1:
   KeysRemaining:=KeysRemaining-1: fi:
NumberTrials:=NumberTrials+1:
od:
return(NumberTrials):
end proc: 
\end{spverbatim}} \normalsize{}
\medskip

The strategy developed in section \ref{Stra2} can be implemented like:

\small{\begin{spverbatim}
Strategy2:=proc(T,N)
local Res,KeyInProgress,LocksRemaining,n,j,Tab:
Tab:=T:  n:=nops(Tab): Res:=0:
KeyInProgress:=1: LocksRemaining:=N:
while LocksRemaining>0 do
if Tab[KeyInProgress]>LocksRemaining then 
   Res:=Res+ LocksRemaining: KeyInProgress:=KeyInProgress+1:
   else  Res:=Res+ Tab[KeyInProgress]: LocksRemaining:=LocksRemaining-1:  
   if KeyInProgress<n then for j from (KeyInProgress+1) to n do 
   if (Tab[j]>Tab[KeyInProgress]) then Tab[j]:=Tab[j]-1: fi: od: fi:
	 KeyInProgress:=KeyInProgress+1: 
fi: od: return(Res):
end proc: 
\end{spverbatim}} \normalsize{}
\medskip

A similar code can be implemented for strategy developed in section \ref{Stra1}. Results of the simulations and comparison with the theoretical results are:

\begin{center}\label{Fig3}
\includegraphics[width=16cm]{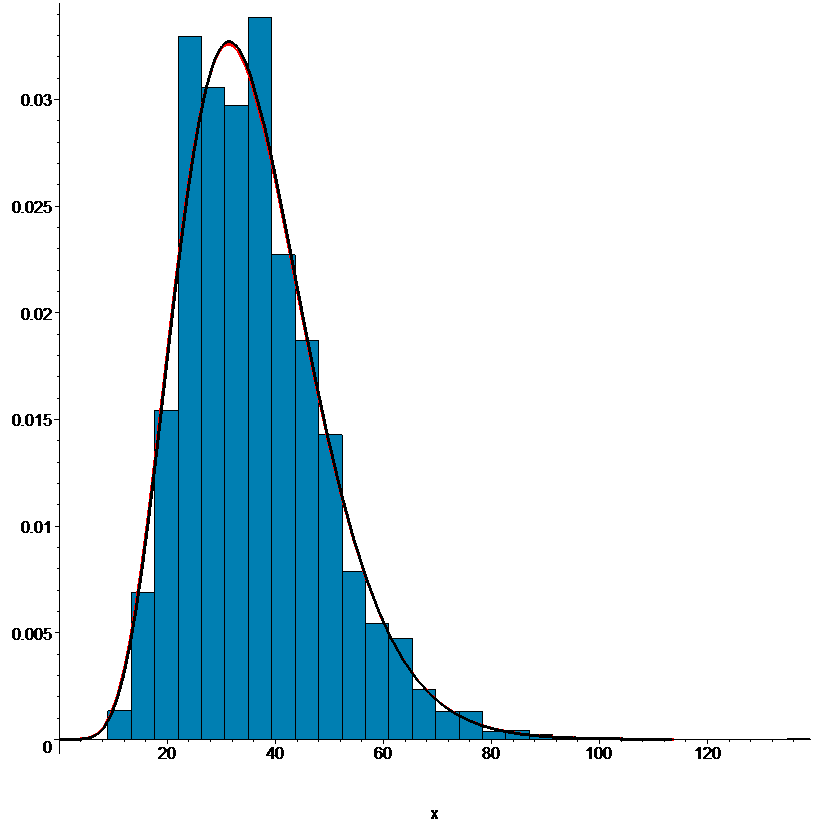}

\textit{Fig. $3$: Illustration of the totally random strategy in the case $n=N=8$. Curve in red represents the $\Gamma$ distribution with the theoretical parameters. Curve in black represents the $\Gamma$ distribution with parameters inferred from the simulation}
\end{center}

\begin{center}\label{Fig4}
\includegraphics[width=16cm]{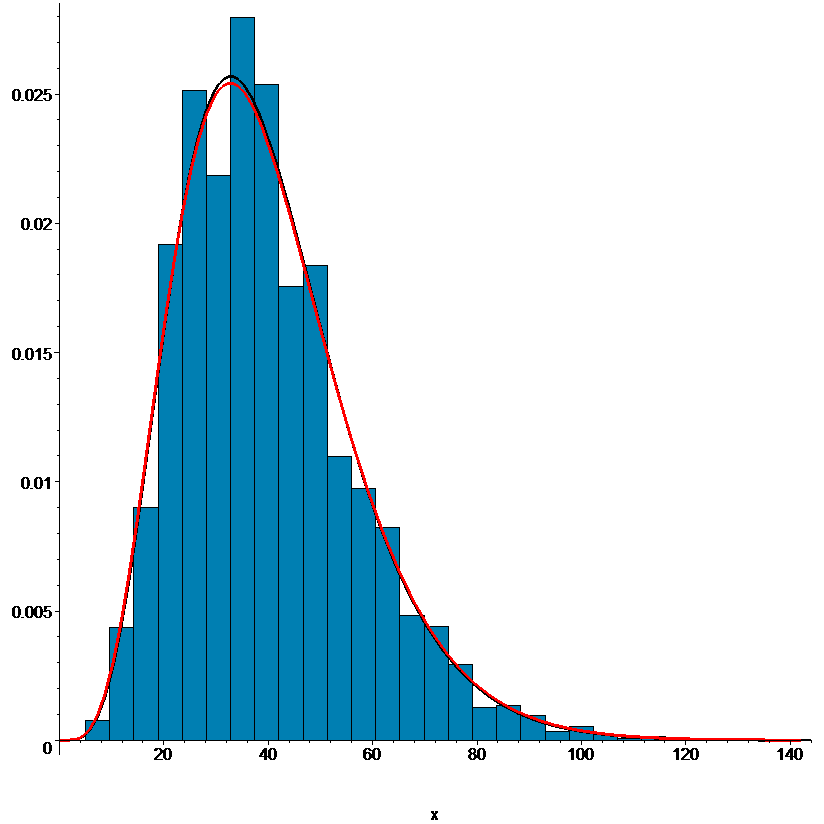}

\textit{Fig. $4$: Illustration of the totally random strategy in the case $n=5$ and $N=10$. Curve in red represents the $\Gamma$ distribution with the theoretical parameters. Curve in black represents the $\Gamma$ distribution with parameters inferred from the simulation}
\end{center}

\begin{center}\label{Fig5}
\includegraphics[width=16cm]{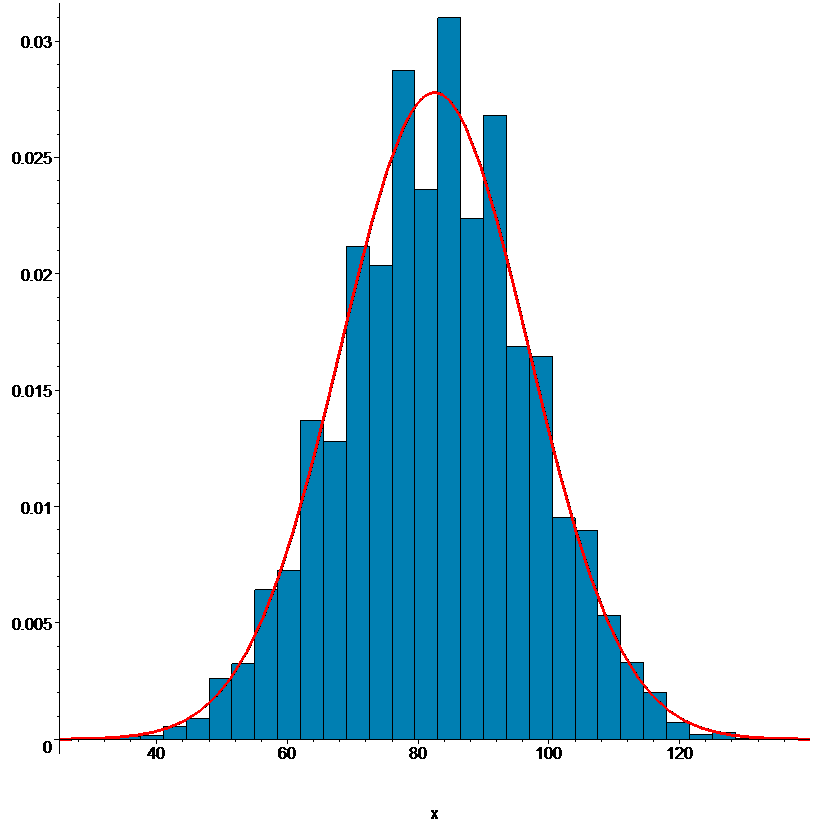}

\textit{Fig. $5$: Illustration of the strategy developed in section \ref{Stra2} in the case $N=20$ and $n=10$. Curve in red represents the normal distribution with the theoretical parameters. Curve in black represents the normal distribution with parameters inferred from the simulation}
\end{center}

\begin{center}\label{Fig6}
\includegraphics[width=16cm]{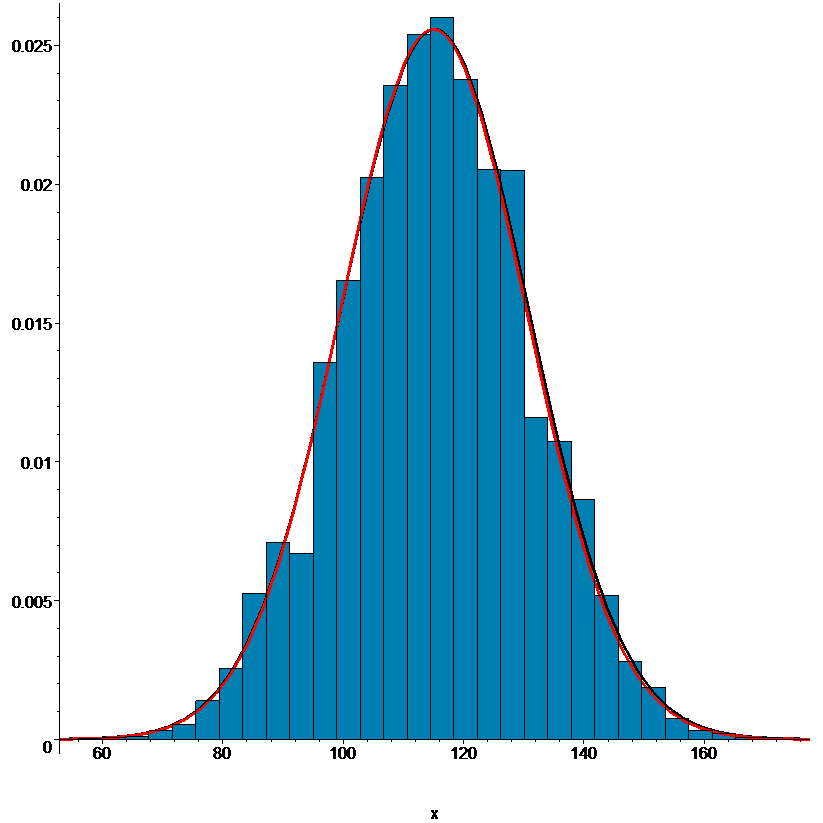}

\textit{Fig. $6$: Illustration of the strategy developed in section \ref{Stra1} in the case $n=N=20$. Curve in red represents the normal distribution with the theoretical parameters. Curve in black represents the normal distribution with parameters inferred from the simulation}
\end{center}

One can obviously observe that the theoretical results are in perfect agreement with the simulations. Contrary to the totally random case, $T_{n,N}$ is better approximated with a classical normal distribution in the other two strategies at least when $n$ is large enough. This is quite logical since it seems reasonable that a Central Limit Theorem may be applied for equation \eqref{T}.  

\section{Conclusion and Outlooks}

Though the problem can be easily stated and understood by common people (after inquiries, this was one of the points that seduced the students to choose this problem among others), its resolution presents various levels of mathematics from elementary ones (small values of $n$ and $N$, totally random strategy) to more involved ones (Chu-Vandermonde identity). The problem is also well suited for simulations and implementation even for students with basic notions of programming. Basic statistical skills are also required to represent the results of simulations. Retrospectively, in addition to the Chu-Vandermonde identity and its connection to the problem, the hardest part for students (even with help) was to unravel the equivalence presented in figure $2$.

\medskip

The problem raises numerous questions and possible developments:  
\begin{itemize}
\item What happens if a key is allowed to open more than one lock? Does it affect the order in which we start over when a lock is opened? What happens if some locks are allowed to be the same?
\item For low $n$ and $N$, it seems possible to get the complete probability distribution of $T_{n,N}$.
\item Simulations in the totally random case indicates that the distribution of $T_{n,N}$ should be close to a $\Gamma$ distribution and arguments where presented to explain this observation. Can these arguments be precised and a complete proof given?
\item When $n$ is large, simulations shows that the distribution of $T_{n,N}$ approaches a normal distribution. Equation \eqref{T} provides a perfect setting to apply some Central Limit Theorem. Can we have a detailed proof of the application of the CLT and a bound relatively to the error committed in the approximation? 
\item Illustration provided by figure $2$ shows that the problem can be recast into a combinatorial problem of counting top-left boxes after selecting some at random (imposing one per line). Hence it would not be surprising that our problem is equivalent to another one in a totally different setting in relation with random permutations.
\item In this article we presented the more natural strategies that come to mind, but are the strategies presented here really the optimal strategies to minimize in average $T_{n,N}$?
\item In real life, changing locks may be more time-consuming than changing keys if the locks are located far apart. Conversely, the keyring might be hard to handle so that changing from a key to the next one may be complicated. Hence one could introduce different times when changing keys or locks. Are the results presented here globally unchanged?
\end{itemize}

\end{document}